\ifcase2 
\documentclass[11pt,a5paper]{article}
\IfFileExists{ajr.sty}{\usepackage{ajr}}
  {\usepackage[a5paper,margin=6mm,tmargin=15mm]{geometry}}
\or 
\documentclass[12pt,a4paper]{article}
  \usepackage[bmargin=20mm,tmargin=23mm]{geometry}
  \AtBeginDocument{%
      \RaisedNamesfalse
      
      }
\or
\documentclass[12pt,a4paper]{article}
  \usepackage[bmargin=20mm,tmargin=23mm]{geometry}
  \AtBeginDocument{%
      \RaisedNamesfalse
      
      }
  \AtEndDocument{\bibliography{bibexport}}
\fi

\title{Accurate and efficient multiscale simulation of a heterogeneous elastic beam via computation on small sparse patches}

\author{%
A.J. Roberts\thanks{School of Mathematical Sciences, University of Adelaide,  South Australia. 
\protect\url{mailto:ProfAJRoberts@protonmail.com}
\protect\url{https://orcid.org/0000-0001-8930-1552}}
\and 
Thien Tran-Duc\thanks{School of Mathematical Sciences, University of Adelaide,  South Australia.
\protect\url{https://orcid.org/0000-0002-2004-5156}}
\and 
J.E. Bunder\thanks{Mathematical Sciences, University of South Australia, Australia.
\protect\url{https://orcid.org/0000-0001-5355-2288}}
\and 
Yannis Kevrekidis\thanks{Departments of Chemical and Biomolecular Engineering \& Applied Mathematics and Statistics, Johns Hopkins University, Baltimore, Maryland, USA.
\protect\url{https://orcid.org/0000-0003-2220-3522}}
}

\usepackage{animate}
\usepackage{microtype,amssymb,amsmath,doi,graphicx,parskip,defns}
\usepackage[defaultlines=3,all]{nowidow}
\usepackage[leftcaption,raggedright]{sidecap}
\usepackage{pgfplots}\pgfplotsset{compat=newest}
\def\tikzsetnextfilename#1{}
\usetikzlibrary{decorations.markings}
\usetikzlibrary{arrows}

\usepackage{natbib}
\usepackage{mycleveref}

\Vec f\Vec n\Vec q\Vec u
\Vec F
\newcounter{i}

\def\symBox#1{#1}
\def\oSym{\symBox{$\color{green!70!black}\circledcirc$}}
\def\xSym{\symBox{$\color{green!70!black}\otimes$}}
\def\uSym{\symBox{$\color{blue}\blacktriangleright$}}
\def\vSym{\symBox{$\color{red}\blacktriangle$}}

\allowdisplaybreaks
\begin{document}

\maketitle


\begin{abstract}
Modern `smart' materials have complex microscale structure, often with unknown macroscale closure.  The Equation-Free Patch Scheme empowers us to non-intrusively, efficiently, and accurately simulate over large scales through computations on only small well-separated patches of the microscale system.  Here the microscale system is a solid beam of random heterogeneous elasticity.  The continuing challenge is to compute the given physics on just the microscale patches, and couple the patches across un-simulated macroscale space, in order to establish efficiency, accuracy, consistency, and stability on the macroscale.  Dynamical systems theory supports the scheme.  This research program is to develop a systematic non-intrusive approach, both computationally and analytically proven, to model and compute accurately macroscale system levels of general complex physical and engineering systems.  
\end{abstract}

\tableofcontents

\section{Introduction}

In structural engineering, microscale lattice materials can be light and highly stiff with customizable macroscale mechanical properties \cite[e.g.,][]{Somnic2022}.  
The challenge we address herein is to accurately and efficiently predict macroscale characteristics emergent from the microscale lattice.
Similarly, composite materials and structures are inherently heterogeneous and anisotropic across multiple scales.
Multiscale modelling is thus critical to the design of composite structures for lightweight mechanical performance \cite[e.g.,][]{Raju2021, Lucarini2022}.
Such composite materials are used in electronics, space, medical, transportation, and other industries \cite[e.g.][]{Matous2017}.
Herein we establish that the Equation-Free Patch Scheme can non-intrusively, efficiently, and accurately simulate over macroscales through computations on only small well-separated patches of the microscale system.

%
%

\begin{figure}[b]
\centering
\caption{\label{figbeamFull}%
movie of a full-domain simulation of a heterogeneous beam showing that beam bending waves and longitudinal compression waves propagate with some `average' properties.  }
\animategraphics[controls,width=\linewidth,every=3]{10}{Figs/beamFullLong/i}{1}{301}
\end{figure}

Consider an example elastic beam with heterogeneous elasticity in 2D as in \cref{figbeamFull}: say 628\,cm long, 20\,cm wide.
The beam is heterogeneous because it is constructed from a modern material with micro-structure of size 3\,cm---so that the heterogeneity is `visible' in \cref{figbeamFull}. 
With a 3\,cm micro-grid, the modelling requires circa 5\,000 variables.
This specific scenario is easily computable, \verb|ode23| took 14\,s \textsc{cpu} time to simulate one period of beam bending oscillation.
But if a more realistic 3\,mm micro-structure is simulated, then the computation time increases by a factor of~1000.  
If 3D elasticity modelling is required for the beam, then the computation time increases by even more orders of magnitude.
The patch scheme \cite[e.g.,][]{Samaey08} we develop herein potentially reduces macroscale computation time by orders of magnitude---more reduction in higher-D space and/or smaller micro-scale.

The patch scheme achieves efficiency by only computing on small sparse patches in space.
\cref{secSnifw} discusses how the patch scheme is non-intrusive in that it just `wraps around' a user's microscale code---a desirable property also identified by \cite{Biezemans2022b}.
The patch scheme, alternatively called the gap-tooth method, ``has formal similarity with \textsc{sp} [superparametrization]'' \cite[p.62]{Majda2014} that was developed in meteorology for weather and climate predictions, and is also akin to the so-called \textsc{fe-fft} and \textsc{fe}\(^2\) methods \cite[e.g.,\S4.7]{Lucarini2022}.

\paragraph{A given microscale discretisation of heterogeneous elasticity}

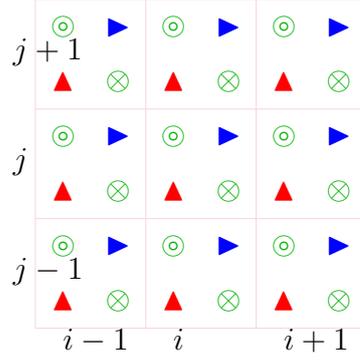
\begin{SCfigure}[1.2]
\centering
\caption{\label{fig:microgrid}a small part of the microscale grid used to code 2D elasticity.  The grid is staggered on the microscale: \uSym, horizontal displacements and velocities; \vSym, vertical displacements and velocities; \oSym, \xSym, components of strain and stress tensor~\eqref{eqs:stress}.
\vspace*{1\baselineskip}
}
\setlength{\unitlength}{4ex}
\def\N{3}\def\Np{4}\def\NN{6}
\def\M{3}\def\Mp{4}\def\MM{6}
\begin{picture}(\NN.5,\MM.5)
\def\symBox#1{\makebox(0,0){#1}}
\put(1,0){
{\color{magenta!20}
    \multiput(-0.5,0.5)(0,2){\Mp}{\line(1,0){\NN}}
    \multiput(-0.5,0.5)(2,0){\Np}{\line(0,1){\MM}} 
    }
\multiput(1,2)(2,0){\N}{\multiput(0,0)(0,2){\M}{\uSym}}
\multiput(0,1)(2,0){\N}{\multiput(0,0)(0,2){\M}{\vSym}}
\multiput(0,2)(2,0){\N}{\multiput(0,0)(0,2){\M}{\oSym}}
\multiput(1,1)(2,0){\N}{\multiput(0,0)(0,2){\M}{\xSym}}
\setcounter{i}0
\multiput(0,0.1)(2,0){\N}{%
    $i\ifcase\value{i}-1\or\or+1\fi$\stepcounter{i}}
\setcounter{i}0
\multiput(-0.9,1.4)(0,2){\M}{%
    $j\ifcase\value{i}-1\or\or+1\fi$\stepcounter{i}}
}
\end{picture}
\end{SCfigure}

We adopt a simple robust microscale approximation of 2D elasticity within the beam.  
On the \emph{staggered} microscale \(xy\)-grid of \cref{fig:microgrid} define the displacements: 
\uSym,~horizontal~\(u_{ij}(t)\); 
\vSym,~vertical~\(v_{ij}(t)\).
Microscale elasticity here first uses centred finite differences to compute stresses, for heterogeneous Lam\'e parameters~\(\lambda,\mu\), at the labelled microscale grid-points (\cref{fig:microgrid}):
\begin{subequations}\label{eqs:stress}%
\begin{align}&
\text{\xSym}\quad \sigma_{xy}
:=\mu_{ij}\big[\delta_ju_{ij}/\delta y_j 
+\delta_iv_{ij}/\delta x_i\big];
\\&
\text{\oSym}\quad \sigma_{xx}
:=(\lambda_{ij}+2\mu_{ij}){\delta_iu_{ij}}/{\delta x_i} 
+\lambda_{ij}{\delta_jv_{ij}}/{\delta y_j};
\\&
\text{\oSym}\quad \sigma_{yy}
:=\lambda_{ij}{\delta_iu_{ij}}/{\delta x_i} 
+(\lambda_{ij}+2\mu_{ij}){\delta_jv_{ij}}/{\delta y_j}.
\end{align}
\end{subequations}
Second, centred finite differences compute the following acceleration \ode{}s
\begin{subequations}\label{eqs:accel}%
\begin{align}
&\text{\uSym}\quad 
\ddot u_{ij}={\delta_i\sigma_{xx}}/{\delta x_i} 
+{\delta_j\sigma_{xy}}/{\delta y_j}\,,
\\&
\text{\vSym}\quad  
\ddot v_{ij}={\delta_i\sigma_{xy}}/{\delta x_i} 
+{\delta_j\sigma_{yy}}/{\delta y_j}\,,
\end{align}
\end{subequations}
potentially with optional small phenomenological damping supplied by a discretisation of~\(\kappa\delsq\dot u_{ij},\,\kappa\delsq\dot v_{ij}\).
The patch scheme wraps around whatever microscale code a user supplies---here it is the microscale system~\eqref{eqs:stress,eqs:accel}

We nondimensionalise the system so that the density is one, and the speed of a macroscale compression wave along the beam is about one, that is, time in these simulations is roughly in milli-seconds.

\paragraph{Random periodic heterogeneity}

\begin{SCfigure}\centering
\caption{\label{fig:heteroE}example of the 2D microscale heterogeneous Young's modulus~\(E_{ij}\) used in computing the elastic Lam\'e parameters~\eqref{eq:Lame}.  In this example, we choose the heterogeneity to have microscale period four along the beam.
\vspace*{2\baselineskip}
}
\includegraphics[scale=0.75]{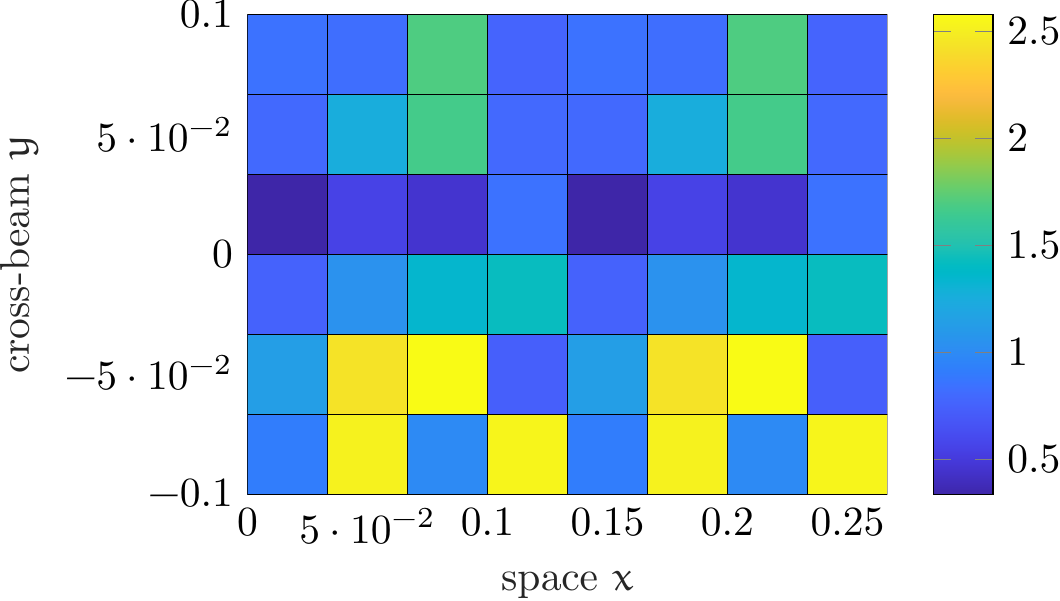}
\end{SCfigure}

The Lam\'e parameters which appear in the stresses~\eqref{eqs:stress} are
\begin{equation}
\lambda:=\frac{\nu E}{(1+\nu)(1-2\nu)},\quad
\mu:=\frac{E}{2(1+\nu)},
\label{eq:Lame}
\end{equation}
in terms of Young's modulus~\(E\) and Poisson ratio~\(\nu\).
To have strong microscale heterogeneity we choose these parameters randomly so that at each microscale grid-point (iid):
\(E_{ij}\)~is log-normal (here varies by factor of about ten); and
\(\nu_{ij}\)~is uniform on~\([0.25,0.35]\).
\cref{fig:heteroE} shows an example~\(E_{ij}\).
Despite such strong heterogeneity, the movie of \cref{figbeamFull} shows the macroscale dynamics appears relatively simple.

\begin{figure}[b]
\centering
\caption{\label{figbeamSim}%
movie of a patch scheme simulation of a heterogeneous beam showing the macroscale propagation across the patches of beam bending waves and longitudinal compression waves.  }
\animategraphics[controls,width=\linewidth,every=3]{10}{Figs/beamSim/i}{1}{301}
\end{figure}

\section{Equation-free patch scheme}

Instead of computing the entire beam as seen in \cref{figbeamFull}, the patch scheme computes only in small sparse spatial patches such as \cref{figbeamSim}.
In this example case, the patch scheme reduces compute time by a factor\({}\propto r:={}\)(patch size)/(spacing~\(H\)), which here is just a modest factor of~\(1/4\).
But with greater scale separation and/or in higher spatial dimensions, the scheme often reduces computational time by many orders of magnitude.

The movie of \cref{figbeamSim} shows a slow progressive wave of beam bending, together with a not-so-slow compression wave along the beam.
These \emph{macroscale predictions are accurate} (\cref{secShpa}) due to the correctness of our simple coupling between patches---even when heterogeneity is strong.
The patch scheme makes these accurate macroscale predictions even when the macroscale closure is unknown:  the scheme does not code a closure.  
Further, `the closure' varies depending upon human assumptions such as choosing averaged models versus cosserat models---the patch scheme makes no such closure assumptions.
The only assumption is that the macroscale quantities of importance vary smoothly between neighbouring patches.

\subsection{Scheme is non-intrusive functional `wrapper'}
\label{secSnifw}

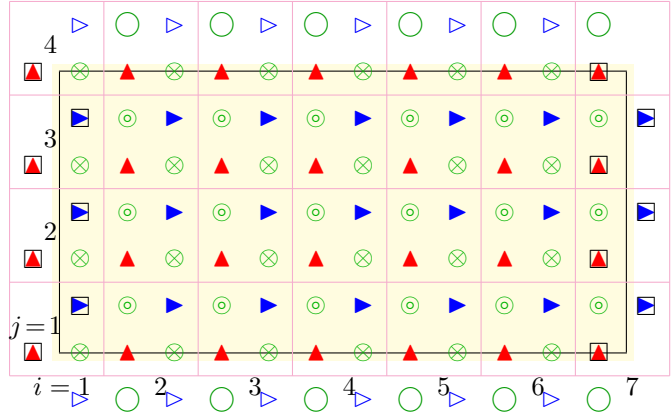
\begin{SCfigure}[1][t]
\centering
\caption{\label{figpatchgridv} 
one example patch of the 2D elastic beam showing the microscale staggered grid (\cref{fig:microgrid}).
This is case of \(n_{\text{subpatch}}=7\) micro-grid intervals along the patch,  and \(n_y=4\) intervals across the beam.
\vspace*{2\baselineskip}
}
{\footnotesize%
\setlength{\unitlength}{4.1ex}%
\def\N{7}\def\Np{8}\def\Nm{6}\def\NN{14}\def\NNmm{12} 
\def\M{4}\def\Mp{5}\def\Mm{3}\def\MM{8} 
\def\symBox#1{\makebox(0,0){#1}}%
\begin{picture}(\NN,\MM.5)
\put(0.5,0){
\put(0.4,1){\colorbox{yellow!15}{\framebox(\NNmm,6){}}}
{\color{magenta!40}
    \multiput(-0.5,0.5)(0,2){\Mp}{\line(1,0){\NN}}
    \multiput(-0.5,0.5)(2,0){\Np}{\line(0,1){\MM}} 
    }
\multiput(1,2)(2,0){\N}{\multiput(0,0)(0,2){\Mm}{\uSym}}
\multiput(0,1)(2,0){\N}{\multiput(0,0)(0,2){\M}{\vSym}}
\multiput(2,2)(2,0){\Nm}{\multiput(0,0)(0,2){\Mm}{\oSym}}
\multiput(1,1)(2,0){\Nm}{\multiput(0,0)(0,2){\M}{\xSym}}
\multiput(0,0)(0,\MM)2{%
  \multiput(1,0)(2,0){\Nm}{\symBox{$\color{blue}\vartriangleright$}}
  \multiput(2,0)(2,0){\Nm}{\symBox{$\color{green!60!black}\bigcirc$}}
  }
\multiput(0,1)(\NNmm,0)2{\multiput(0,0)(0,2){\M}{\symBox{$\square$}}}
\multiput(1,2)(\NNmm,0)2{\multiput(0,0)(0,2){\Mm}{\symBox{$\square$}}}
\setcounter{i}0
\multiput(0,0.1)(2,0){\N}{\stepcounter{i}%
    $\ifnum\value{i}=1i=\else\quad\fi\arabic{i}$}
\setcounter{i}0
\multiput(-0.5,1.4)(0,2){\M}{\stepcounter{i}%
    $\ifnum\value{i}=1j\!=\!\else\quad\;\fi\arabic{i}$}
}
\end{picture}
}
\end{SCfigure}

Consider one of the patches of the 2D beam shown in \cref{figbeamSim}.
With the \emph{given} microscale \(xy\)-grid (\cref{fig:microgrid}), zooming in to the microscale each patch is like that of \cref{figpatchgridv}.
Here each patch extends across the cross-section (\(y\)-dimension) of the beam. 
Open symbols in \cref{figpatchgridv} are ghost nodes outside the
patch and implement \emph{given} stress-free top\slash bottom conditions on the beam.
The \emph{only} addition required by the patch scheme are the edge values (`squared' micro-grid nodes in \cref{figpatchgridv}) on the left\slash right of each patch.

The patch scheme couples patches together by providing the patch-edge values through interpolation across the macroscale between patches \cite[e.g.,][]{Roberts06d, Roberts2011a, Cao2014a}.  
Here we interpolate from each of the \emph{centre patch values across} the beam (\(i=4\) in \cref{figpatchgridv}) of `nearby' patches, to determine the corresponding patch-edge value.
Here we implement spectral (\textsc{fft}) interpolation between the patches for high accuracy (\cref{secShpa}).
The scheme does \emph{not} presume that any average is appropriate.

This implementation shows that the patch scheme is non-intrusive \cite[e.g.,][]{Biezemans2022b}: it just `wraps around' any micro-grid code a user trusts.  
Consequently, we provide a toolbox \cite[]{Maclean2020a} for others to implement the patch scheme around their micro-code.

\subsection{Scheme embeds macroscale dynamics}

Given the patch scheme does not assume anything about what are `correct' macroscale variables, \emph{a crucial question} is the following: how can we be assured that the patch scheme captures the macroscale slow dynamics?
\quad\emph{An answer} is provided by the \cite{Whitney1936} embedding theorem.

Roughly, the theorem is that every \(m\)D~manifold is parametrisable from almost every subspace of more than~\(2m\)D. 
Let's see what this means for us.
In essence, the patch scheme provides the higher-D subspace in which the slow manifold of the macroscale wave dynamics is embedded.  

For beams in two spatial dimensions, the basic macroscale beam models have, at each cross-section,  displacement and velocity of both bending and compression.  
Thus the elastic beam dynamics has a slow manifold that is \(m=4\)D at every cross-section.\footnote{Such statements, invoking a manifold or subspace ``at every cross-section", are in a sense developed by the theory of \cite{Roberts2013a}.  That is, in systems of large spatial extent there often are important, spatially global, invariant manifolds of high-D that are effectively decomposable into a union of \emph{spatially local} manifolds\slash subspaces of relatively lower dimension---a dimension determined by the spatial cross-section---and that are weakly coupled to neighbouring locales.}
Alternatively, 2D cosserat beam models add a shear mode to the macroscale model---two more variables---leading to a not-quite-so-slow manifold of \(m=6\)D at every cross-section.
These physically based models are \emph{slow manifolds} because they focus on the relatively slow waves of solutions varying slowly in space, and neglect all the faster high-frequency cross-waves.

In the patch scheme, \cref{figbeamFull,figbeamSim} show simulations with a cross-section of \(n_y=7\) micro-grid intervals, but let's discuss the case of just \(n_y=4\) (\cref{figpatchgridv}).
For \(n_y=4\), there are seven microscale nodes across each patch edge.  Each node has a displacement and velocity, and so leads to a \(14\)D~subspace for macroscale communication between patches.

Because \(14>2\cdot6>2\cdot4\)\,, the Whitney embedding theorem asserts that the patch scheme exchanges enough information to almost surely parametrise both such slow manifolds of the macroscale dynamics.
The patch scheme does \emph{not} need to explicitly compute and exchange  specific assumed macroscale average quantities.

%
%
%

\section{Scheme has proven accuracy}
\label{secShpa}

\cref{secmapc} discusses established theory which generally proves that the patch scheme makes accurate macroscale predictions.
Such proofs are in stark contrast to the vast \emph{machine learning\slash artificial intelligence} developments which prove very few general results: for example, \cite{Brenner2021} comment ``\ldots\textsc{ml} studies, as the lack of rigorous theory does not offer (yet!) guarantees of convergence''.
Before discussing theory, we first report some computational verification \text{of high accuracy.}

\subsection{Computation verifies exactness}

\begin{SCfigure}
\centering
\caption{\label{fig:1vis}multiscale spectrum of eigenvalues~\(\lambda\) separates macroscale modes on the right from sub-patch microscale modes on the left.  The axes are scaled nonlinearly.
Here the small viscosity is~\(0.001\) so the microscale decays, but the macroscale waves are long-lasting.
\vspace*{4\baselineskip}
}
\input{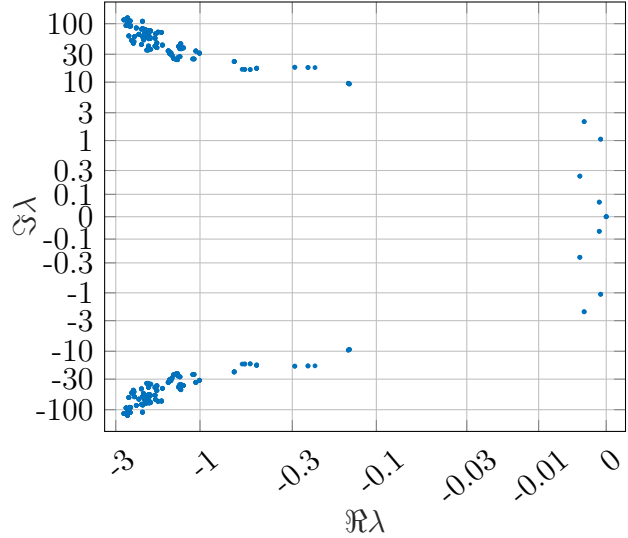}
\end{SCfigure}
Here we restricted attention to linear elasticity so we know that the wrapped patch system is fully characterised by the resultant Jacobian matrix. 
We numerically compute the Jacobian matrix of the patch scheme by elementary numerical differentiation.

Because of the macroscale translational invariance of the patch scheme, the macroscale eigenvectors are correctly sinusoidal.  Hence the only macroscale errors occur in the eigenvalues of the Jacobian. 
\cref{fig:1vis} plots the spectrum of all eigenvalues for one example of random heterogeneity, in the case of five patches for simplicity.
Observe there are:
\begin{itemize}
\item (on the right) four \(\lambda=0\) of rigid beam motion;
\item four \(-0.001 \pm\i1.057\) and four \(-0.003 \pm\i2.111\) of compressions waves;
\item four \(-0.001\pm\i0.061\) and four \(-0.004\pm\i0.237\) of beam bending waves;
\item with the above \emph{macroscale} eigenvalues separated by a \emph{spectral gap} from the following sub-patch microscale eigenvalues;
\item (on the left) many \(\Re\lambda<-0.1\) of {uninteresting} sub-patch micro-scale fast-waves (headed by ten eigenvalues around \(-0.14\pm\i9.29\)).
\end{itemize}
To quantify the accuracy, \cref{tblerr} compares eigenvalues obtained from full-domain code, with the above macroscale eigenvalues obtained by the wrapped patch scheme.
For all patch size ratios and heterogeneities tested, the patch scheme's macroscale eigenvalues are exact to numerical round-off error.

\begin{SCtable}
\centering
\caption{\label{tblerr}error in patch scheme's macroscale eigenvalues~\(\lambda\) for various patch size ratios~\(r\):  the macroscale~\(\lambda\)s are \emph{exact} to round-off error---due to patch coupling by \emph{spectral} interpolation.  
}
\begin{tabular}{llll}
\hline
macro-eigenvalue & \(r=\tfrac12\) & \(r=\tfrac14\) & \(r=\tfrac18\) 
\\[.5ex]\hline
\(-0.001\pm\i0.061\) & 2e-12 & 1e-12 & 2e-13\\
\(-0.001\pm\i0.061\) & 2e-12 & 4e-12 & 2e-12\\
\(-0.004\pm\i0.237\) & 1e-12 & 8e-13 & 3e-12\\
\(-0.004\pm\i0.237\) & 1e-12 & 2e-12 & 3e-12\\
\(-0.001 \pm\i1.057\) & 7e-13 & 4e-13 & 6e-13\\
\(-0.001 \pm\i1.057\) & 6e-13 & 5e-13 & 6e-13\\
\(-0.003 \pm\i2.111\) & 1e-13 & 2e-13 & 2e-13\\
\(-0.003 \pm\i2.111\) & 4e-13 & 5e-13 & 2e-13\\\hline
\end{tabular}
\end{SCtable}

Such exactness is due to the spectral interpolation used here.
If, instead of spectral, local polynomial interpolation of degree~\(p\) is used to couple the patches, then generally the patch scheme has macroscale errors\({}\propto H^{p}\) where \(H={}\)inter-patch spacing \cite[e.g.,][]{Roberts06d, Roberts2011a}.

\begin{SCfigure}
\centering
\caption{\label{fig:0vis}multiscale spectrum of eigenvalues~\(\lambda\) for the patch scheme in the case of zero viscosity.  
The horizontal axis shows that all modes have zero real-part to numerical round-off error.
That is, in the case of zero viscosity, \emph{this patch scheme preserves the wave nature of the underlying physics.}
\vspace*{4\baselineskip}
}
\input{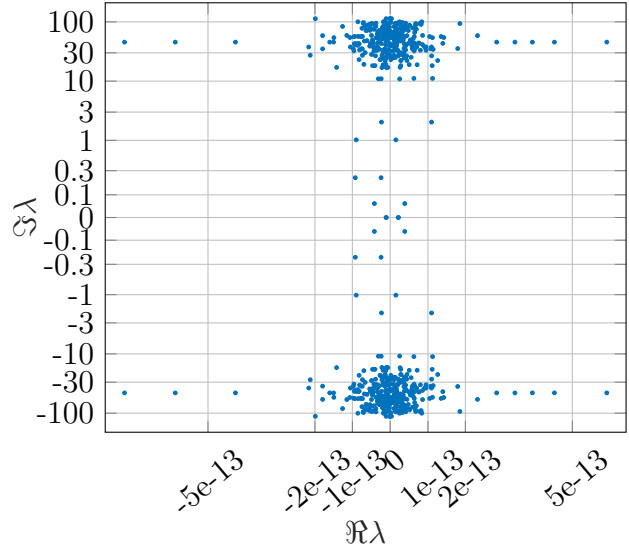}
\end{SCfigure}

\paragraph{Undamped waves?}
With zero viscosity, there are only oscillations in the underlying physics. 
In such a scenario computational methods are very delicate.  
Here, \cref{fig:0vis} illustrates that all eigenvalues of the patch scheme have \(|\Re\lambda|<10^{-12}\).\footnote{In some realisations of the heterogeneity, the sensitive multiplicity four eigenvalue \(\lambda=0\) numerically splits into four showing~\(|\Re\lambda|\) up to~\(10^{-6}\) due to round-off errors.}
Hence, even with no viscosity, the \emph{patch scheme preserves the oscillatory wave nature of the heterogeneous physics.}

There is a perception that the patch scheme ``only works well on problems with an inertial manifold and for systems in which most modes are strongly decaying'' \cite[p.62]{Majda2014}.  
This verification of accuracy for purely elastic beams shows that this perception is false.
Applications and theory for other wave systems also refute this perception \cite[e.g.,][]{Cao2014a, Bunder2020a, Divahar2022b}.

\subsection{Mathematical analysis proves consistency}
\label{secmapc}

Mathematical analysis has proven properties of the patch scheme in general. Mostly, the published proofs explicitly address dissipative (nonlinear) systems.
However, as discussed by \cite{Bunder2020a}, the patch scheme in space only recasts spatial interactions, so whether the time derivative is~\(\D t{}\) of dissipation or~\(\DD t{}\) of waves makes little difference.

Two complementary types of results have been proven. 
They involve the spacing between patch centres~\(H\).
First, Centre Manifold Theory may be applied at finite spacing~\(H\) by introducing a `bookkeeping' parameter~\(\gamma\) to label inter-patch communication \cite[e.g.,][\S2]{Roberts2011a} to prove the existence of a slow manifold in the patch scheme (including when it is applied to nonlinear systems).
Then the parameter~\(\gamma\) structures inter-patch interactions, and their algebraic expression, to empower theory based at \(\gamma=0\), via \emph{regular} perturbation, to address finite~\(\gamma\) such as the case of full coupling \(\gamma=1\) \cite[e.g.,][Cor.~2]{Roberts2011a}.

Second, the patch scheme is \emph{consistent} with the underlying micro-code as the patch spacing \(H\to0\) \cite[e.g.,][Thm.~7]{Roberts2011a}.
The consistency is that the macroscale of the patch scheme is the same as the macroscale of the given micro-coded system, to errors~\Ord{H^p} when using polynomial interpolation of degree~\(p\).
For example, spectral interpolation corresponds to `\(p=\infty\)' so then errors vanish to all orders as in \cref{tblerr}.


These results and general proofs were first done for homogeneous systems \cite[e.g.,][]{Roberts06d, Roberts2011a}.
They were subsequently extended to heterogeneous microscales \cite[]{Bunder2013b}, and recently extended to alternative inter-patch coupling that preserves self-adjointness \cite[]{Bunder2020a}.
Interestingly, the extension of the theoretical support to heterogeneous cases invokes the ensemble of all phase-shifts of the heterogeneity.
The ensemble is spatially homogeneous, so the homogeneous proofs and results apply to establish the heterogeneous results.

%

\section{Conclusion}

As an initial exploration of the patch scheme for homogenisation of heterogeneous elasticity, we considered the prototypical case of a 2D elastic beam.  
The scheme gives a non-intrusive and efficient \emph{computational homogenisation} of \emph{given} microscale system via spatially sparse small patches.
The patch coupling has proven accuracy, controllable error, at \emph{finite} scale separation.

The patch scheme makes only one assumption: in the scenarios of interest to a user, there is no significant spatial structures in the mesoscale between the patch spacing~\(H\) and the microscale resolved in the patches.
In contrast to most other multiscale methods, there is:  no assumed boundary conditions on Representative Volume Elements (variously periodic, Dirichlet, Neumann);  no explicitly assuming slow variables;  and no presumed necessary variational principle.
The scheme is entirely \emph{physically interpretable}: there is no hidden mystic machinations of neural networks \cite[e.g.,][]{Brenner2021}

The patch scheme is simple to apply.
In contrast to other multiscale methods there is: no arbitrary averaging; no oversampling regions; no buffer regions; no action regions; no guessed fast/slow variables; no epsilons; and no limits. 
As a non-intrusive `wrapper', anyone can start using the patch scheme via a \textsc{Matlab}\slash Octave Toolbox \cite[]{Maclean2020a, Roberts2019b}

\paragraph{Acknowledgements} 
This research was supported by Australian Research Council grants DP220103156 and DP200103097.

\end{document}